\documentclass[12pt]{article}

\usepackage{a4}   
\usepackage{amssymb,amsmath,mathrsfs,textcomp}   
\usepackage[amsmath,thmmarks]{ntheorem}   
\usepackage{yhmath}   
\usepackage{indentfirst}   
\usepackage{tocloft}   
\usepackage[hyperfootnotes=false,colorlinks=true,linkcolor=black,citecolor=black]{hyperref}   

\usepackage{graphicx}
\usepackage{epstopdf}

\theoremstyle{nonumberplain}
\theorembodyfont{\rm}
\theoremseparator{.}
\newtheorem{definition}{Definition}

\theoremstyle{plain}
\theorembodyfont{\rm}
\theoremseparator{.}

\theoremstyle{plain}
\theorembodyfont{\slshape}
\theoremseparator{.}

\theoremstyle{plain}
\theorembodyfont{\slshape}
\theoremseparator{$^\prime$.}

\theoremstyle{nonumberplain}
\theorembodyfont{\slshape}
\theoremseparator{.}

\theoremstyle{nonumberplain}
\theoremheaderfont{\bf}
\theorembodyfont{\rm} 
\theoremseparator{.}
\theoremsymbol{$\Box$}
\newtheorem{proof}{Proof}

\makeatletter
\renewcommand*{\@seccntformat}[1]{\csname the#1\endcsname.\quad}
\makeatother

\newcommand{\coloneqq}{\mathrel{\mathop:}=}
\newcommand{\eqqcolon}{=\mathrel{\mathop:}}

\makeatletter
\def\moverlay{\mathpalette\mov@rlay}
\def\mov@rlay#1#2{\leavevmode\vtop{%
   \baselineskip\z@skip \lineskiplimit-\maxdimen
   \ialign{\hfil$\m@th#1##$\hfil\cr#2\crcr}}}
\newcommand{\charfusion}[3][\mathord]{
    #1{\ifx#1\mathop\vphantom{#2}\fi
        \mathpalette\mov@rlay{#2\cr#3}
      }
    \ifx#1\mathop\expandafter\displaylimits\fi}
\makeatother

\begin{document}

%
%
\title{A domain with non-plurisubharmonic squeezing function}
\author{John Erik Forn{\ae}ss$^{(\ast)}$ and Nikolay Shcherbina}
\maketitle
\nopagebreak

%
%
\small\noindent{\bf Abstract.}
We construct a strictly pseudoconvex domain with smooth boundary whose squeezing function is not  plurisubharmonic. \normalsize

%
%
\renewcommand{\thefootnote}{}\footnote{$^{(\ast)}$The first author was supported in part by the Norwegian Research Council grant number 240569 and NSF grant DMS1006294}

%
%

%
%
%
%
%
\section{Introduction}

In this paper we are dealing with the properties of squeezing functions on domains. The idea of using this concept goes back to the papers  \cite{LiuSunYau04} and  \cite{LiuSunYau05} where a new notion of holomorphic homogeneous regular domains was introduced. The last kind of domains can be seen as a generalization of Teichm{\"u}ller spaces, and, as it was shown in \cite{LiuSunYau04}, \cite{LiuSunYau05} and \cite{Yeung09}, they admit many nice geometric and analytic properties.

Motivated by the mentioned above works \cite{LiuSunYau04} and  \cite{LiuSunYau05},  Deng, Guan and Zhang in \cite{DengGuanZhang12} introduced the notion of squeezing functions defined for arbitrary bounded domains:

\begin{definition}
Let $\Omega$ be a bounded domain in $\mathbb C^n$. For $p \in \Omega$ and a
holomorphic embedding $f: \Omega \rightarrow {\mathbb B}^n$ satisfying $f(p) = 0$ we set
\[ S{_{\Omega}}(p, f) \coloneqq \rm{sup}\{r > 0 : r {\mathbb B}^n \subset f(\Omega)\}, \]
and then we set
\[ S{_{\Omega}}(p) \coloneqq \rm{sup}_{f} \{S{_{\Omega}}(p, f)\}, \]
where the supremum is taken over all holomorphic embeddings  $f: {\Omega} \rightarrow {\mathbb B}^n$ with $f(p) = 0$ and ${\mathbb B}^n$ is representing the unit ball in $\mathbb C^n$. The function $S{_{\Omega}}$ is called the squeezing function of ${\Omega}$.
\end{definition}

Properties of the squeezing function for different classes of domains were then studied in \cite{DengGuanZhang12}, \cite{DengGuanZhang16} and \cite{KimZhang15}. Moreover, using the results of \cite{DiederichFornaessWold14}, sharp estimates not only for the squeezing functions, but also for the Carath\'{e}odory,
Sibony and Azukawa metrics near the boundary of a given strictly pseudoconvex domain were obtained in \cite{FornaessWold15}. Similar results for the Bergman metric are given in \cite{DiederichFornaess15}.

On the other hand, in many cases functions which are naturally defined on pseudoconvex domains enjoy plurisubharmonicity properties (see, for example, \cite{Yamaguchi89} and \cite{Berndtsson06}). That is why a few years ago the following question was raised: 

{\it Is it always true that the squeezing function of a strictly pseudoconvex domain with smooth boundary is plurisubharmonic?}

The main result of this paper gives a negative answer to the question and can be formulated as follows.

\noindent \textbf{Theorem.} {\it There exists a bounded strictly pseudoconvex domain with smooth boundary in $\mathbb C^2$ whose squeezing function is not plurisubharmonic}.

\section{Preliminaries} \label{sec_preliminaries}

First we briefly recall the definitions of the Kobayashi and Carath\'{e}odory metrics.
Let $\Delta$ denote the unit disc, and let ${\mathcal{O}}(M,N)$ denote the set of holomorphic maps
from M to N. For a domain $\Omega \subset \mathbb C^n$ we consider an arbitrary point $p \in \Omega$ and an arbitrary vector $\xi \in \mathbb C^n$.

$\bullet$ {\it Kobayashi metric} $K{_{\Omega}}(p, \xi)$. We define

$K{_{\Omega}}(p, \xi) = \rm{inf}\{|{\alpha}|  ; {\,\, \exists} \,\, f \in {\mathcal{O}}(\Delta, \Omega)  \,\,\,\,  f(0) = p, \,\, {\alpha} f'(0) = \xi\}.$

$\bullet$ {\it Carath$\acute{e}$odory metric} $C{_{\Omega}}(p, \xi)$. We define

$C{_{\Omega}}(p, \xi) = \rm{sup}\{|f'(p) (\xi)|; {\,\, \exists} \,\, f \in {\mathcal{O}}(\Omega, \Delta) \,\,\,\,    f(p) = 0\}.$

\vspace{2mm}
\noindent 
Observe that the above definitions imply directly the next well known properties of metrics.

\vspace{3mm}
\noindent \textbf{Monotonicity of Metrics.} {\it Let $\Omega_1 \subset \Omega_2$ be bounded domains in $\mathbb C^n$, $p$ be a point in $\Omega_1$ and $\xi$ be an arbitrary vector in $\mathbb C^n$. Then the following properties hold true}
\[ K{_{\Omega_1}}(p, \xi) \geq  K{_{\Omega_2}}(p, \xi)\,\,\,\,\,\,\,  and  \,\,\,\,\,\,\,\,\,\,\,\,C{_{\Omega_1}}(p, \xi) \geq  C{_{\Omega_2}}(p, \xi). \]
We will also need the following two statements which one easily gets from the definitions (detailed proofs of them can be found in \cite{DengGuanZhang12}).

\vspace{3mm}
\noindent \textbf{Lemma 1.} {\it Let $\Omega$ be a bounded domain in $\mathbb C^n$. Then for all $z \in \Omega$ and all $\xi \in \mathbb C^n$ one has}
\[ S{_{\Omega}}(p) K{_{\Omega}}(p, \xi) \le C{_{\Omega}}(p, \xi) \le K{_{\Omega}}(p, \xi). \]
\noindent \textbf{Lemma 2.} {\it The squeezing function $S{_{\Omega}}$ of any bounded domain $\Omega$ in $\mathbb C^n$ is continuous.}

\vspace{3mm}
The last statement implies, in particular, the following property (a slightly weaker result was stated as Theorem 2.1 in \cite{DengGuanZhang16}, but a slight modification of the proof presented there gives actually the stronger statement as it is formulated below).

\vspace{3mm}
\noindent \textbf{Lemma 3.} {\it Let $\Omega$ be a bounded domain in $\mathbb C^n$. Then for any compact set $K \subset \Omega$ and any $\epsilon > 0$ there exists $\delta > 0$ such that for each subdomain $\widetilde{\Omega}$ of $\Omega$,   $K \subset \widetilde{\Omega}$, having the property that $b{\widetilde{\Omega}} \subset U_{\delta}(b{\Omega})$ one has $|S{_{\Omega}}(p) - S{_{\widetilde{\Omega}}}(p)| < \epsilon$ for every $p \in K$. Here by $U_{\delta}(b{\Omega})$ is denoted the $\delta$-neighbourhood of the boundary $b{\Omega}$ of $\Omega$.}

\vspace{3mm}

Now we give some estimates on the Carath\'{e}odory and Kobayashi metrics of some special domains.

\vspace{3mm}
\noindent \textbf{Lemma 4.} {\it Let $\,\,\,\, 0 < a  < 1 < b < + \infty$ be given numbers. For each $m \in {\mathbb N}$, consider the domain
\[ \Omega'_{m} \coloneqq \{(z, w) \in {\mathbb C^2}: a < |z| < b, |w| < 1, |w| <  |z|^{-{m}} \}. \]
Then there exists $C >0$ such that $C_{\Omega'_m}(p, \xi) \le C$ for $p = (1, 0)$, $\xi = (1, 1)$ and all $m \in {\mathbb N}$. }

\begin{figure}[!htb]
\centering
\includegraphics[trim= 0cm 0.8cm 0cm 2.4cm, clip=true, scale=1.0]{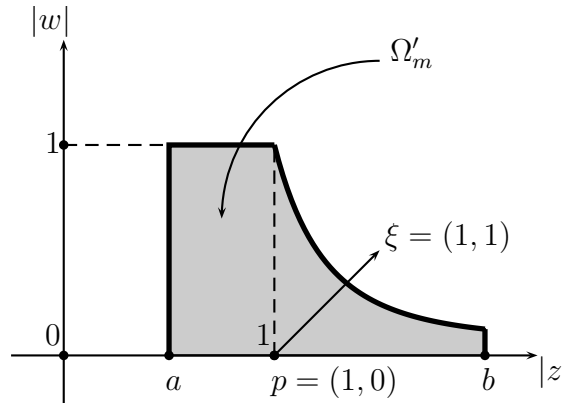}
\caption{The domain $\Omega_m'$.}
\label{fig1}
\end{figure}

\begin{proof}
Consider an arbitrary function $f \in {\mathcal{O}}(\Omega'_{m}, \Delta)$ such that $f(p) = 0$. Observe that the restriction $f_v$ of $f$ to the vertical disc $\Delta_v \coloneqq \{z = 1\} \times \{|w| < 1\} = \{z = 1\} \cap {\Omega'_m}$ centered at $p$ is a holomorphic function from $\Delta_v$ to $\Delta$ having the property $f_v(p) = 0$. Then, by the Schwarz lemma, one has 
\[ |f'(p)(0, 1)| = |{\frac{\partial f}{\partial w}}(p)| = |f'_v(p)| \leq 1. \]
Similarly, for the restriction $f_h$ of $f$ to the horizontal disc 
\[ \Delta_h \coloneqq \{|z - 1| < \rm{min}(1 - a, b - 1)\} \times \{w = 0\} \subset {\Omega'_m} \cap \{w = 0\} \]
we have that $f_h :  \Delta_h \rightarrow \Delta$ is a holomorphic function such that $f_h(p) = 0$. Hence, in view of the Schwarz lemma, one also has
\[ |f'(p)(1, 0)| = |\frac{\partial f}{\partial z}(p)| = |f'_h(p)| \leq \frac{1}{\rm{min}(1 - a, b - 1)}. \]
Therefore
\[ |f'(p) (1, 1)| = |\frac{\partial f}{\partial z}(p) + \frac{\partial f}{\partial w}(p)| \leq |\frac{\partial f}{\partial z}(p)| + |\frac{\partial f}{\partial w}(p)| \leq  \frac{1}{\rm{min}(1 - a, b - 1)} + 1 \eqqcolon C. \]
Since $f$ was an arbitrary function from ${\mathcal{O}}(\Omega'_m, \Delta)$ such that $f(p) = 0$, we finally conclude that for the Carath\'{e}odory metric the estimate $C{_{\Omega'_m}}(p, \xi) \le C$ holds true for all $m \in {\mathbb N}$.
\end{proof}

\vspace{-0.3cm}
\noindent \textbf{Lemma 5.} {\it For each $m \in {\mathbb N}$, consider the domain
\[ \Omega''_m \coloneqq \{(z, w) \in {\mathbb C^2}: |w| < 1, |w| < |z|^{-m} \}. \]
Then $K{_{\Omega''_m}}(p, \xi) \geq \sqrt{\frac{m}{2}}$ for $p = (1, 0)$, $\xi = (1, 1)$ and each $m \in {\mathbb N}$. }

\begin{figure}[!htb]
\centering
\includegraphics[trim= 0cm 0.8cm 0cm 2.6cm, clip=true, scale=0.9]{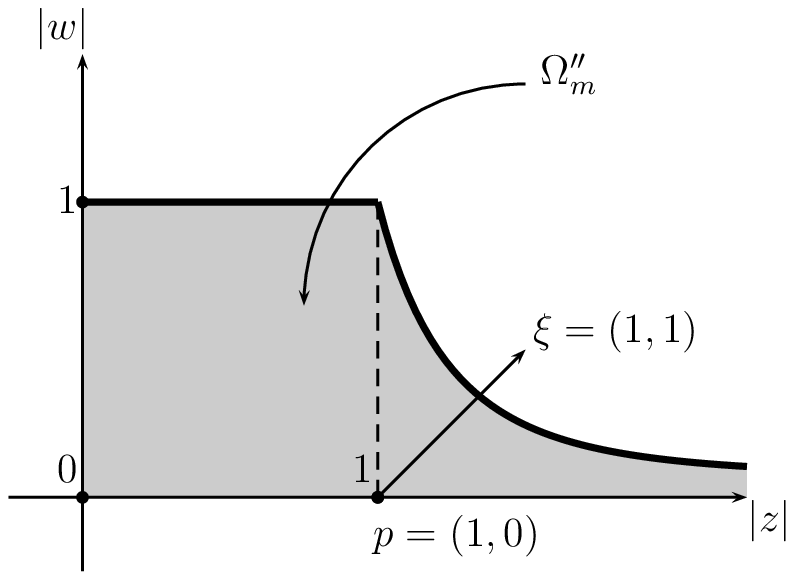}
\caption{The domain $\Omega_m''$.}
\label{fig2}
\end{figure}

\begin{proof}
Consider an arbitrary map $f \in {\mathcal{O}}(\Delta, \Omega''_{m})$ such that $f(0) = p$ and ${\alpha} f'(0) = \xi =(1, 1)$ for some $\alpha$. Then $f$ can be represented by 
\[ f(\zeta) = (z(\zeta), w(\zeta)) = (1  + \frac{1}{\alpha} \zeta + a_2 {\zeta}^2 +  ...,  \frac{1}{\alpha} \zeta + b_2 {\zeta}^2 +  ...), \] 
where $\zeta \in \Delta$. Since, by the definition of $\Omega''_m$, one has $|w z^m| < 1$, it follows that 
\[ 1 > | (\frac{1}{\alpha} \zeta + b_2 {\zeta}^2 +  ... ) ( 1  + \frac{1}{\alpha} \zeta + a_2 {\zeta}^2 +  ... )^m | = | \frac{1}{\alpha} \zeta + (b_2 + \frac{m}{{\alpha}^2}) {\zeta}^2 +  ... |  . \]
Then, from the Schwarz type bound for higher order coefficients (see Theorem 2 in \cite{Ruscheweyh85} for a relatively recent generalization of the classical Schwarz inequality to similar bounds for all coefficients of the Taylor expansion), we get that
\begin{equation} \label{schwarz1}
|b_2 + \frac{m}{{\alpha}^2}| \leq 1.
\end{equation}
Since, by the definition of $\Omega''_m$, one also has
\[ |\frac{1}{\alpha} \zeta + b_2 {\zeta}^2 +  ...| = |w| \leq 1, \]
we conclude from the mentioned above Schwarz type bound for the higher order coefficients that
\begin{equation} \label{schwarz2}
|b_2| \leq 1.
\end{equation}
Combining estimates (1) and (2), we get
\[ |\frac{m}{{\alpha}^2}| \leq 2 \,\,\,\, \Rightarrow \,\,\,\, |\alpha| \geq \sqrt{\frac{m}{2}}, \]
\noindent
which gives the desired estimate $K{_{\Omega''_m}}(p, \xi) \geq \sqrt{\frac{m}{2}}$ for each $m \in {\mathbb N}$.
\end{proof}

\section{Example} \label{sec_example}

We first construct an auxiliary domain which we will denote by $\Omega$. Let $a > 1$ be an arbitrary number, which will be fixed in what follows, and let $1 < a_1 < a_2 < ... < a_k < ... < a$ be a sequence (which will also be fixed) such that $\lim_{k \to \infty}{a_k} = a$. We define $\Omega$ as the set of points $(z, w) \in \{\frac{1}{a}<|z|< a\} \times {\mathbb{C}}_w$ satisfying the following conditions:
\[ \left\{ \begin{array}{c@{\,,\quad}l} |w| < {B_k} \, |z|^{n_k}  &  {\rm{for}} \,\,  \frac{1}{a_{k + 1}} < |z| \leq \frac{1}{a_k}, \, k = 1, 2, 3, ... , \\ |w| < 1 &   {\rm{for}} \,\,  \frac{1}{a_1} < |z| \leq a_1,   \\ |w| < {B_k} \, |z|^{-{n_k}} &   {\rm{for}} \,\,  a_k  \leq |z| < a_{k + 1}, \, k = 1, 2, 3, ...   \end{array} \right.\]

The numbers $n_k$ and $B_k$ will be defined inductively so that $B_1 = 1$, and for each $k \in \mathbb{N}, k \geq 2,$ one has $n_{k} > n_{k - 1}$ and ${B_{k - 1}} \, {a_k}^{-n_{k - 1}} = {B_{k}} \, {a_k}^{-n_{k}}$ (the last condition guarantees that the functions defining $\Omega$ will match at the points $a_k$ and $\frac{1}{a_k}$, $k \in \mathbb{N}$) and, moreover, the inequality $S_{\Omega}(p_k) < \frac{1}{k}$ for the squeezing function on $\Omega$ at the point $p_k = (a_k, 0)$ holds true for every $k \in \mathbb{N}$.\\

\begin{figure}[h]
\centering
\includegraphics[trim= 0cm 0.9cm 0cm 2.6cm, clip=true, scale=1.0]{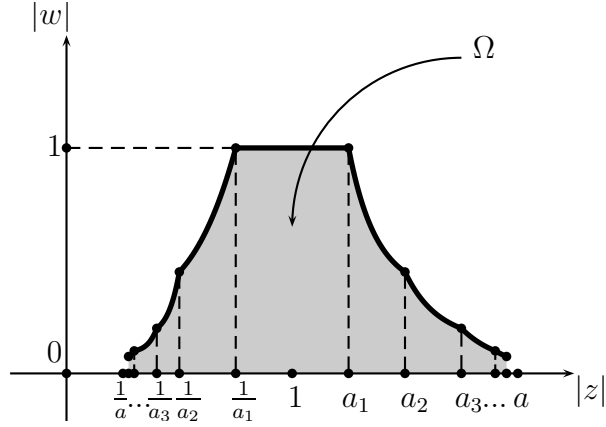}
\caption{The auxiliary domain $\Omega$.}
\label{fig3}
\end{figure}

The starting point of our inductive construction is the definition of $\Omega$ over the annulus $\{ \frac{1}{a_1} < |z| \leq a_1 \}$ by the inequality $|w| < 1$. Now we describe the inductive step of this construction. Assume that the part $\Omega_k$ of the domain $\Omega$ over the annulus $\{\frac{1}{a_k}<|z|< a_k\}$ is already constructed, i.e., we have already defined the numbers $n_q, B_q$ for $q = 1, 2, ..., k - 1$. For being able to find suitable values of $n_k$ and $B_k$, we first make a biholomorphic change of coordinates $F_k$ in ${\mathbb{C}}^{\ast} \times \mathbb{C}$:
\[ z \rightarrow \frac{z}{a_k} \eqqcolon z', \, w \rightarrow w \frac{{a_{k}}^{n_{k - 1}}}{{B_{k - 1}}}  \Big({\frac{z}{a_k}}\Big)^{n_{k - 1}} \eqqcolon w'. \]
Observe that in new coordinates $(z', w')$ the part of the domain $F_k(\Omega)$ over the annulus $\{ \frac{a_{k - 1}}{a_k} \leq |z'| < 1 \}$ is defined by $|w'| < 1$ and the part of $F_k(\Omega)$ over the annulus $\{ 1 \leq |z'| < \frac{a_{k + 1}}{a_k} \}$ is defined by $|w'| < |z'|^{-({n_{k}} - {n_{k - 1}})}$, where  $n_{k}$ still has to be chosen. Note  also that the domain  
\[ F_k(\Omega \cap (\{ a_{k - 1}  < |z| < a_{k + 1} \} \times {\mathbb{C}}_w) ) = \]
\[ = \{(z', w'):\frac{a_{k - 1}}{a_k} < |z'| < \frac{a_{k + 1}}{a_k},|w'| < 1,|w'| < |z'|^{-({n_k} - {n_{k - 1}})} \}  \]
has the form $\Omega'_m$ (see Lemma 4 for the description of $\Omega'_m$) with $m = {n_k} - {n_{k - 1}}, a = \frac{a_{k - 1}}{a_k}, b = \frac{a_{k + 1}}{a_k}$ and it is a proper subdomain of the domain $F_k(\Omega)$. Moreover, since for each $k \in \mathbb{N}$ the inequality $n_{k} > n_{k - 1}$ holds, the domain $F_k(\Omega)$ will be contained in the domain $\Omega''_m$ (see Lemma 5 for the description of $\Omega''_m$) with $m = n_{k} - n_{k - 1}$. 

\begin{figure}[!htb]
\centering
\includegraphics[trim= 0cm 0.8cm 0cm 2.6cm, clip=true, scale=1.0]{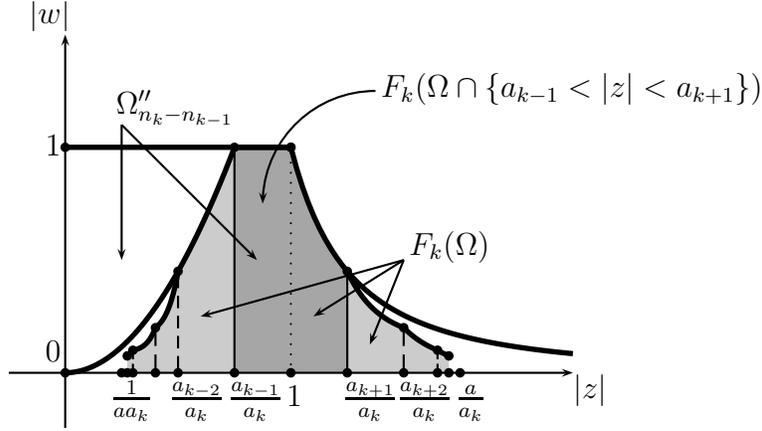}
\caption{The domains $F_k\big(\Omega \cap \{a_k < |z| < a_{k+1}\} \big)$, $F_k(\Omega)$ and $\Omega_{n_k-n_{k-1}}''$.}
\label{fig4}
\end{figure}

Hence, in view of monotonicity of the Carath\'{e}odory metric and Lemma 4, one has
\[ C{_{F_k(\Omega)}}(p, \xi) \le C{_{\Omega'_{m}}}(p, \xi) \le C_k \]
for $p = (1, 0)$, $\xi = (1, 1)$ and all $m \in {\mathbb N}$. We also have from monotonicity of the Kobayashi metric and Lemma 5 that
\[ K{_{F_k(\Omega)}}(p, \xi) \geq K{_{\Omega''_{m}}}(p, \xi) \geq \sqrt{\frac{m}{2}} \]
for $p = (1, 0)$, $\xi = (1, 1)$ and each $m \in {\mathbb N}$. It follows then from Lemma 1 that
\[ S{_{F_k(\Omega)}}(p) \le \frac{C{_{F_k(\Omega)}}(p, \xi)}{K{_{F_k(\Omega)}}(p, \xi)} \le C_k \sqrt{\frac{2}{m}} \]
and hence $S{_{F_k(\Omega)}}(p) < \frac{1}{k}$ for $n_{k} > n_{k - 1} + 2 {k^2}{C_k^2}$. If we choose now $n_{k}$ satisfying the last inequality, then, using the condition $ {B_{k - 1}} \, {a_k}^{-{n_{k - 1}}} = {B_k} \, {a_k}^{-{n_k}}$, we can easily compute ${B_{k}} = {B_{k - 1}} \, {a_k}^{{n_{k}} - {n_{k - 1}}}$. Finally, note that, in view of biholomorphic invariance of the squeezing function,
\[ S{_{\Omega}}(a_k) =  S{_{F_k(\Omega)}}(p) < \frac{1}{k}, \]
for each $k \in \mathbb{N}$. This completes the inductive step of our construction of the auxiliary domain $\Omega$.

Now we are ready to construct a strictly pseudoconvex domain with non-plurisub-harmonic squeezing function. Note first that $\Omega$ is pseudoconvex by construction. Observe also that, since the map $z \rightarrow \frac{1}{z}, w \rightarrow w$ is a biholomorphic automorphism of $\Omega$, and, since the squeezing function is biholomorphically invariant, one has
\[ S{_{\Omega}}\Big(\frac{1}{a_k}\Big) =  S{_{\Omega}}(a_k) < \frac{1}{k}, \]
for each $k \in \mathbb{N}$. Take now $p = (1, 0) \in \Omega$, denote $c  \coloneqq S{_{\Omega}}(p) >0$ and fix from now on a number $k \in \mathbb{N}$ so large that $\frac{1}{k} < c$. Then, using Lemma 3 with $\epsilon < \frac{1}{2} \, (c - \frac{1}{k})$, we approximate the domain $\Omega$ from inside by a strictly pseudoconvex smoothly bounded domain $\widetilde{\Omega}$ (one can obviously choose this domain to be also circular in $z$ and $w$) so well that for every point $q$ of the set 
\[ (\{|z| = \frac{1}{a_k} \} \times \{w = 0\}) \cup (\{|z| = a_k \} \times \{w = 0\}) \subset \widetilde{\Omega} \cap \{w = 0\} \]
one has
\[ S{_{\widetilde{\Omega}}}(q) < \frac{1}{k} + \epsilon < c - \epsilon < S{_{\widetilde{\Omega}}}(p). \]
This means that the maximum principle for the restriction of the function $S_{\widetilde{\Omega}}(\bf \cdot)$ to the annulus $\{{\frac{1}{a_k}} \le |z| \le  {a_k}\} \times \{w = 0\} \subset {\widetilde{\Omega}} \cap \{w = 0\}$ does not hold and, hence, the function $S_{\widetilde{\Omega}}(\bf \cdot)$ cannot be plurisubharmonic. Thus ${\widetilde{\Omega}}$ is a strictly pseudoconvex domain as desired. The proof of the Theorem is now completed. $\hfill \Box$

\vspace{3mm}
\noindent
{\bf Remark.} In the proof above instead of using Lemma 3 it is enough to  use the weaker statement of   Theorem 2.1 from \cite{DengGuanZhang16} at the points ${\frac{1}{a_k}}$, $a_k$ and $p$ and the circular invariance of the domain ${\widetilde{\Omega}}$ and the squeezing function $S_{\widetilde{\Omega}}(\bf \cdot)$.

\vspace{3mm}
\noindent
{\bf Acknowledgement.} {\it Part of this work was done while the second author was a visitor at the
Capital Normal University (Beijing). It is his pleasure to thank this institution
for its hospitality and good working conditions.}

%
%
 \vspace{1truecm}

%
%
%
{\sc J. E. Forn{\ae}ss: Department of Mathematics, NTNU --- 7491 Trondheim, Norway}
  
{\em e-mail address}: {\texttt john.fornass@math.ntnu.no}\vspace{0.3cm}
  
{\sc N. Shcherbina: Department of Mathematics, University of Wuppertal --- 42119 Wuppertal, Germany}
  
{\em e-mail address}: {\texttt shcherbina@math.uni-wuppertal.de}

\end{document}